\documentclass[12pt]{article}
\usepackage{amsmath,amssymb,amsthm,array,amsfonts}
\usepackage{mathrsfs} 
\usepackage[english]{babel}
\usepackage[dvipdfm,bookmarks,colorlinks]{hyperref}
\usepackage[dvips,final]{graphicx}

\theoremstyle{definition}
\newtheorem{dfn}{Definition}
\newtheorem{defn}{Definition}
\theoremstyle{theorem}
\newtheorem{thm}{Theorem}

\newtheorem{cor}{Corollary}
\newtheorem{prp}{Proposition}
\newtheorem{rem}{Remark}
\newtheorem{prop}{Proposition}

\begin{document}
\begin{titlepage}
  \begin{center}
    \textbf{National Taras Shevchenko University of Kyiv}
  \end{center}\vspace{2.5cm}
\begin{center}
    {\bf \Large Set-membership state estimation framework for uncertain linear
      differential-algebraic equations}  
\end{center}\vspace{1.5cm}
\begin{center}
    {\bf Sergey Zhuk}
\end{center}\vspace{1.5cm}
\begin{center}
Technical report 06BF015-02/2008
\end{center}\vspace{1.5cm}
\begin{center}
  Kyiv-2008
\end{center}
\end{titlepage}

%\maketitle
% \section{ToDo}
% \begin{itemize}
% \item check why pics are not shown (compare to automatica paper)  
% \item spelling
% \end{itemize}
\begin{center}
    {\bf Set-membership state estimation framework for uncertain linear
      differential-algebraic equations}  
\end{center}\vspace{.5cm}
\begin{center}
    {\bf Sergey Zhuk}
\end{center}\vspace{1.5cm}
\textbf{Abstract.}  We investigate a state
estimation problem for the dynamical system described by uncertain linear operator
equation in Hilbert space. The uncertainty is supposed to admit a
set-membership description. We present explicit expressions for
linear minimax estimation and error provided that any pair of uncertain
parameters belongs to the quadratic bounding set. We introduce a new notion of
minimax directional observability and index of non-causality for linear
noncausal DAEs. Application of these notions to the state
estimation problem for linear uncertain noncausal DAEs allows to derive new minimax
recursive estimator for both continuous and discrete time. We illustrate the benefits of
non-causality of the plant applying our approach to scalar nonlinear
set-membership state estimation problem. Numerical example is presented.  


\section{Introduction and problem statement}
The applications of differential-algebraic equations
(DAEs or descriptor systems) in economics, demography, mechanics and engineering
are well known~\cite{Lewis1986}. This in turns motivates researchers to
investigate DAEs from the mathematical point of view~\cite{Mehrmann2005}. Here
we focus on a design of state estimation 
algorithm for uncertain linear non-causal DAE.  

The most common approach to DAEs investigation is
to reduce it to some canonical form which in turn is equal to some normal
ODE. In particular, one of the basic results of the
algebraic theory of regular linear DAEs with constant matrices\footnote{Here
  we present a formulation for the finite-dimensional Banach space.} was
introduced in \cite{Rutkas1975}: if the linear DAE with constant matrices
\begin{equation}
  \label{eq:ldae}
  F\dot x=Cx+Bf
\end{equation}
is well defined ($\mathrm{det}[\lambda F-C]\not\equiv0$ ) then for all initial values
$x(t_0)=x_0$ there exists the unique solution $x(\cdot)$ provided that
$f(\cdot)$ is sufficiently smooth. The index $s$ of the pencil $F,C$
is said to be an index of linear DAE \eqref{eq:ldae}. One can reduce
\eqref{eq:ldae} to the ODE via change of coordinates so that the pencil $F,C$
brings into canonical form \cite{Dai1989} and differentiating exactly $s$
times provided that $f$ is sufficiently smooth. In such a way one can derive
an analogue of the celebrated Cauchy formula for the linear regular DAEs with
constant matrices. This result is generalized to variable coefficients by
means of a standard canonical form (SCF): in \cite{Campbell1983} it was shown
that \eqref{eq:ldae} with analytical $F,C,B$ is solvable (i.e. for every sufficiently smooth $f$ there exists at least one continuously
differentiable solution to \eqref{eq:ldae} provided $F,C$ to be
sufficiently smooth) if there exists the SCF for \eqref{eq:ldae}. 
%  generalizes this result to variable coefficients: provided $F,C$ to be
% sufficiently smooth DAE \eqref{eq:ldae} is said to be solvable at $[t_0,T]$ if
% for every sufficiently smooth $f$ there exists at least one continuously
% differentiable solution. In \cite{Campbell1983} it was shown that
% \eqref{eq:ldae} with analytical $F,C,B$ 
% is solvable iff there exists a standard canonical form for
% \eqref{eq:ldae}.
Note that in this case $\mathrm{rank}\,F(t)$ changes only
 at finite number of points from within any compact $[t_0,T]$. In
 \cite{Campbell1987} it was noted 
that not all solvable DAEs can be put into SCF and the solvable DAE is equal
to some differential-algebraic equation in the canonical form which generalize
SCF. In this respect we say that DAE is causal if it can be reduced -- at least
locally in nonlinear case -- into normal ODE. The geometry of
the reduction procedure for nonlinear causal DAEs $F(x,\dot x)=0$ was
investigated in \cite{Reich1990,Rabier1994}, where the index of DAE was
defined as a smallest natural $s$ so that the sequence of the constraint
manifolds \cite{Reich1990} \begin{equation*}
  \begin{split}
    & M_k:=TW_{k-1}\cap M_{k-1}, M_0:=\{(x,p):F(x,p)=0\},\\
    &W_0:=\{x\in\mathbb R^n:(x,p)\in M_0\}
\end{split}
\end{equation*}
becomes stationary for $k>s$. This coincides with the
definition of the index of linear DAE. Further discussion of the DAEs
solvability theory and related topics is presented in
\cite{Lewis1986,Mehrmann2005}. 

The noncausal DAE differs radically from the causal one. For instance,
consider 
\begin{equation}
  \label{eq:ldae1}
\dfrac d{dt}\bigl[
\begin{smallmatrix}
  1&&0
\end{smallmatrix}\bigr](
\begin{smallmatrix}
  x_1\\x_2
\end{smallmatrix})=\bigl[
\begin{smallmatrix}
  c_1&&c_2
\end{smallmatrix}\bigr](
\begin{smallmatrix}
  x_1\\x_2
\end{smallmatrix})+f(t), \bigl[
\begin{smallmatrix}
  1&&0
\end{smallmatrix}\bigr](
\begin{smallmatrix}
  x_1\\x_2
\end{smallmatrix})(t_0)=x_0
\end{equation} Let $x_2(\cdot)\in\mathbb
L_2(t_0,T)$,  $f(\cdot)\in\mathbb L_1(t_0,T)$ and
$x_0\in\mathbb R$. By definition put $$
x_1(t):=\exp(c_1(t-t_0))x_0+\int_{t_0}^t\exp(c_1(t-s))c_2x_2(s)+f(s)ds
$$ % Clearly, $t\mapsto(x_1(t),x_2(t))^T$ represent a solution for any
% $f(\cdot)\in\mathbb L_1(t_0,T)$ and $x_0\in\mathbb R$.
It is clear that any solution of~\eqref{eq:ldae1} is given by the formula $
t\mapsto(x_1(t),x_2(t))^T$. According to a behavioral approach
\cite{Ilchmann2005} one can think about $x_2$ as an input or as a part of the
system state representing uncertain inner disturbance generated by the plant
itself. % This example shows that non-causal DAE, in general, has non-unique
% solution. Also it is easy to give an example when non-causal DAE may not be
% solvable for any $f(\cdot)$ and $x_0$.
In order to clarify this ambiguity we shall give an exact definition of the
DAEs solution accepted in this paper. According to~\cite{Zhuk2007} $x(\cdot)$
is said to be a solution of  
\begin{equation}
  \label{eq:FDAE}
  \dfrac d{dt}Fx(t)=C(t)x(t)+f(t)
\end{equation} 
with initial condition $Fx(t_0)=0$ if $Fx(\cdot)$ is totally continuous
function, $x(\cdot)$ satisfies \eqref{eq:FDAE} almost everywhere and
$Fx(t_0)=0$ holds. This definition allows to properly define the adjoint
system. Also it guarantees that the linear mapping induced by \eqref{eq:FDAE}
 is closed in corresponding Hilbert space~\cite{Zhuk2007}. Note that this
 doesn't hold for \eqref{eq:ldae}. Another useful application of the introduced
 solution is in the control theory. In \cite{Ozcaldiran1989} authors discuss
 difficulties arising while applying of proportional feedback $f=Kx$ to the
 \eqref{eq:ldae}: even 
well defined DAE ($\mathrm det(sF-C)\ne0$) may become singular ($\mathrm
det(sF-C-BK)\equiv 0$). In \cite{Kurina2007} a properly
stated leading term $A(t)\dfrac d{dt}F(t)x$ is used in order to give a
feedback solution to LQ-control problem with DAE constraints. This generalizes
the definition of DAE solution \cite{Zhuk2007} to the case of variable
matrices.  

Recently solvability conditions for abstract semi-linear non-causal DAE has been
studied in \cite{Rutkas2008} assuming that the pencil $sF-C$ is singular,
$F,C$ are closed linear mappings in abstract Banach space. Properties of the
solutions of noncausal implicit differential equation with special structure
were discussed in \cite{Rabier1989}. Note that
non-causal DAEs are not just a "pure" mathematical structure which is
suitable for solving control or observation problems only -- some potential
applications of non-causal DAEs was briefly discussed in \cite{Grimm1988}. 

A state estimation framework for
linear dynamic models has several widely-used approaches: $H_2/H_\infty$
filtering and set-membership state estimation. $H_2$-estimators like Kalman or
Wiener filters (also known as minimum variance
filters~\cite{Balakrishnan1984}) give estimations of the system state with
minimum error variance. These filters require an exact model of signal
generating process and full information about a statistical nature of noise
sources. Recently, the $H_2$-estimation for linear DAEs has been studied in
\cite{nikh1}. Authors derive a so-called "3-block" form for the optimal
filter and a corresponding 3-block Riccati equation using the maximum likelihood
approach. The obtained recursion is stated in terms of a block matrix
pseudoinverse. In \cite{ishixara2} the filter recursion is represented in
terms of a deterministic data fitting problem solution. Authors introduce an
explicit form of the 3-block matrix 
pseudoinverse for a descriptor model with special structure, so that the form of
obtained in \cite{ishixara2} filter coincides with presented in \cite{nikh1}. 
% $H_2$-estimation for the
% linear time-invariant (LTV) stochastic descriptor systems with regular matrix
% pencil was considered in~\cite{gerdin}.
A brief overview of steady-state $H_2$-estimators is presented
in~\cite{Deng1999}. Optimal $H_\infty$ estimators minimize the $2$-induced
norm\footnote{$\sup$ of the relation  
between the Euclidean norms of estimation errors and model disturbances} of
the operator that maps unknown disturbances with finite energy to filtered
errors~\cite{Basar1995}. In literature it is common to construct
suboptimal estimators~\cite{Shaked1992} that guarantee aforementioned norm to
be less then a prescribed performance level $\gamma$.    
 % This requires to check if given level $\gamma$  is greater or equal then the minimal performance level which is not so easy to compute.
Note that $H_\infty$ estimators are certain Krein space $H_2$
filters~\cite{Sayed1996}. Krein space approach was used in~\cite{Zhang2003}
for risk-sensitive filtering in linear time-invariant (LTI) descriptor models
with regular matrix pencil under stochastic noise. A linear matrix inequality
approach was used in~\cite{Xu2007} in order to construct reduced order
$H_\infty$-filter for LTI DAE with regular matrix pencil. An up to date description of the state of the art is to be found at\cite{Xu2006}. 

% One of the basic notions of the set-membership state estimation framework is
% an a posterioiri set or informational set. By definition it is a set of all
% possible state vectors consistent with measured output data provided that
% input disturbances and measurement errors are some elements of the given
% bounded sets.
% 
% Set-membership filters are based on the notion of the informational set (a
% posteriori set). By definition it is a set of all possible state vectors
% consistent with measured output data when input disturbances and measurement
% errors are some elements of the given bounded sets.
In the sequel %  we focus on
% that case when the state estimation is assigned to the central point of a set
% of possible states. There are several interconnected techniques aimed to
% construct such estimation:
% optimization~\cite{Bertsekas1971,Tempo1985,Nakonechnii1978}, 
% set-valued analysis~\cite{Kurzhanski1997,Bakan2003} and game theory
% approach~\cite{Kuntsevich1992}. For further references see
% \cite{Chernousko1994,Milanese1991}.\\
% Here 
 we focus on the following problem: given
some element (for instance measurements of the system output) $y$
 from some functional space one needs to estimate the expression
 $\ell(\theta)$ provided that $g(\theta)=0$. This
 problem becomes non-trivial if the latter equation has more than one
 solution and the equality $y=C(\theta)$ holds. In this case the
 estimation problem may be reformulated as follows: given
 $y=C(\theta),\theta\in\Theta,y\in Y$ one needs to find the estimation
 $\widehat{\ell(\theta)}$ of the expression $\ell(\theta)$ provided that
 $g(\theta)=0$ and $C(\cdot),\ell(\cdot)$ are given functions. Note that
 $\widehat{\ell(\theta)}:=\ell(\hat\theta)$ if the equation $y=C(\theta)$ has
 the unique solution $\hat\theta$. 
  
The estimation problem is said to be linear if $\Theta,Y$ are linear spaces
and $C(\cdot),\ell(\cdot)$ are linear mappings. It is a common case when $$
C(\theta)=H\varphi+D\eta,
g(\theta)=L\varphi+Bf, \eqno(*)
$$ where $\theta=(x,f,\eta)\subset X\times F\times Y$, $H,D,L,B$ are linear
mappings. The linear estimation problem is said to be uncertain if $D\ne0$,
$L$ and $B$ are non-trivial or if $B=0$ and
$\mathrm{N}(L)=\{\varphi:L\varphi=0\}\ne\{0\}$. Note that the choice of the
solution method depends on the ``type of uncertainty'': if $f,\eta$ denote
realizations of random elements then it's natural to apply probability
methods. This requires an a priori knowledge of distribution characteristics
of the random elements. In the sequel we assume that there is uncertainty in
$(*)$ if distributions of random elements or some deterministic parameters of
the system are partially unknown. It is natural to choose the estimation from
some class to be optimal in the sense of the given criteria. According to this
the linear uncertain estimation problem is said to be minimax if the class of
estimations $\widehat{\ell(\varphi)}$ is restricted to all linear functions
$(u,y)+c$ of $y$ and the criteria is set to be the minimum of the worst-case
error. A description of the state of the art in the theory of linear uncertain
minimax estimation problems with special $\ell,L,H,B,D$ in special spaces is
to be found at \cite{Bertsekas1971,Nakonechnii1978,Tempo1985,Nakonechnii2004,Chernousko1994,Milanese1991,Kurzhanski1997,Bakan2003,Kuntsevich1992}. 
\subsubsection{Author's own research activities} 
Classical theory of uncertain estimation problems \cite{Bertsekas1971}-\cite{Nakonechnii2004} works well when the linear mapping
$L$ in $(*)$ has bounded inverse. One of the author's theoretical achievements is
the extension of the linear minimax estimation theory 
%\cite{Chernousko1994,Milanese1991,Kurzhanski1997,Bakan2003,Kuntsevich1992} 
 to abstract equations with closed linear non-injective mapping \cite{Zhuk2008e,Zhuk2008b,Zhuk2008d,Zhuk2009} in Hilbert space. This extension is based on the general duality principle asserting that the linear minimax
estimation problem is equal to some control problem with convex non-smooth cost
and linear constraints provided that uncertain parameters belong to closed
bounded convex sets in corresponding Hilbert spaces 
\cite{Zhuk2006c, Zhuk2006b,Zhuk2008b}. Note,
that these results were previously obtained for the finite dimensional space 
\cite{Zhuk2004,Zhuk2004a,Zhuk2003,Zhuk2004b,Zhuk2004b,Zhuk2004c}. 

In order
to apply the abstract theory to DAEs the sufficient conditions on DAEs
matrices were introduced asserting that the linear mapping induced by the
noncausal DAE is closed and has closed range \cite{Zhuk2007}. Also a
generalization of the integration by parts formula and the
necessary and sufficient conditions of solvability of DAE in the
form~\eqref{eq:FDAE} is presented in \cite{Zhuk2007}. The solvability
condition is obtained via application of Tikhonov regularization approach. 
With help of these results new notions of the minimax
directional observability and index of causality for discrete time
\cite{Zhuk2009a} and continuous time linear non-causal DAEs
\cite{Zhuk2009b} were introduced. % If system state $x(s)$ is observable in the
% minimax sense in the direction $\ell$ then the projection of the reachability
% set (consistent with measurements $y(t),t_0\ge t\le s$) onto direction $\ell$
% in $\mathbb R^n$ is given by $[-\hat\sigma,\hat\sigma]$, where $\hat\sigma$
% denotes the minimax estimation error. If $x(s)$ is unobservable in
% the minimax sense for $\ell$ then the minimax estimation is set to zero and
% $\hat\sigma=+\infty$. This means that the structure of the measurements do not
% provide any information about $(\ell,x(s))$.
The minimax directional observability provides a qualitative description of
DAEs singularity with a respect to the given observations. Using this notion
author developed several representations of the minimax estimation 
% in terms of the linear 2-point boundary value problem solutions 
\cite{Zhuk2005c,Zhuk2005a,Zhuk2006a}. The final result is an
algorithm which allows to compute the minimax estimation of the linear
noncausal DAE state in the real time \cite{Zhuk2009b}. The structure of this
algorithm coincides with celebrated Kalman filter recursions for normal linear
ODEs with continuous time. Similar results were obtained for linear noncausal
DAEs with discrete time
\cite{Zhuk2008a,Zhuk2006,Zhuk2008c,Zhuk2005b,Zhuk2005,Zhuk2008,Zhuk2009a}. In
\cite{Zhuk2009a} the author gives a complete solution to the problem of
recursive implementation of the minimax a-posteriori estimation (similar to
posed in~\cite{Tempo1985}) for the linear non-causal DAEs with discrete time. 

Also the theory developed in \cite{Zhuk2008b} was applied to the state
estimation of the solutions of finite-dimensional linear boundary-value
problems \cite{Zhuk2007a} and the minimax mean-square estimations of trends
\cite{Zhuk2008f}.

\emph{Notation}. $c(G,\cdot)=\sup\{(z,f),f\in G\}$,\\
$\delta(G,x)=0$ if $x\in G$ and $+\infty$ otherwise,\\
$\mathrm{dom}f=\{x\in\mathcal H:f(x)<\infty\}$,\\
$f^*(x^*)=\sup_{x}\{(x^*,x)-f(x)\}$, $(L^*c)(u)=\inf\{c(G,z),L^*z=u\}$, \\
$(fL)(x)=f(Lx)$, $(L^*c)(u)=\inf\{c(G,z),L^*z=u\}$,\\
$\mathrm{cl}f=f^{**}$, \\
$\mathrm{Arginf}_uf(u)$ denotes the set of minimum points of $f$,\\
$\partial f(x)$ denotes the sub-differential of $f$ at $x$,\\
$(\cdot,\cdot)$ denotes the inner product in Hilbert space,\\
%$\ell_2$ denotes a space of all
%vector-sequences $\{f_k\}=(f_1,f_2,\dots)$ with finite norm $\|{f_k}\|^2=\sum
%(f_k,f_k)$, 
$S>0$ means $(Sx,x)>0$ for all $x$,\\ 
% from within appropriate Hilbert space,
$L^*$ denotes adjoint operator,\\
$P_{L^*}$ denotes the orthogonal projector onto $R(L^*)$,\\
$R(L)$, $N(L)$ and $\mathscr{D}(L)$ denote the range, the null-space and the
domain of the linear mapping $L$,\\
$F'$ denotes transposed matrix,\\
$F^+$ denotes pseudoinverse matrix,\\
$E$ denotes the identity matrix,\\
$\mathrm{diag}(A_1\dots A_n)$ denotes diagonal matrix with $A_i$,
$i=\overline{1,n}$ on its diagonal, \\
$\overline{G}$ denotes the closure of the set $G$,\\
$[x_1,\dots,x_n]$ denotes an element of the Cartesian
product $H_1\times\dots \times H_n$ of Hilbert spaces
$H_i$,$i=\overline{1,n}$,\\
$\mathbb R^n$ denotes $n$-dimensional arithmetic Hilbert space,\\
$\mathrm C^{m\times n}(t_0,T)$ denotes the space of all continuous on
$(t_0,T)$ functions with values in $\mathbb R^{m\times n}$,\\
$\mathbb L_2(t_0,T)$ denotes the space of all measurable functions with finite
integral $\int_{t_0}^Tf^2dt$,\\
$\mathbb W_2^m(t_0,T)$ denotes the space of all absolutely continuous functions with
derivative from $\mathbb L_2(t_0,T)$,\\ 
$M\xi$ denotes expected value of the random vector $\xi$. 
\section{Linear uncertain estimation problem}
In this section we present the main result of the paper \cite{Zhuk2008b}. All
proofs are given in \cite{Zhuk2008b}. 

Suppose that $L\varphi\in\mathscr{G}$ and 
\begin{equation}
  \label{eq:y}
  y=H\varphi+\eta
\end{equation}
The mappings $L,H$ and the set $\mathscr G$ are supposed to be given. The
element $\eta$ is uncertain. Our aim is to solve the inverse problem: to
construct the operator mapping the given $y$ into the estimation
$\widehat{\ell(\varphi)}$ of expression $\ell(\varphi)$ and to calculate the
estimation error $\sigma$. Now let us introduce some definitions. 

The operator $L:\mathcal{H}\mapsto\mathcal{F} $ is assumed to be closed. Its
domain $\mathscr{D}(L)$ is supposed to be a dense subset of the Hilbert space
$\mathcal{H}$, $H\in\mathscr{L}(\mathcal{H},\mathcal{Y})$. Note that the
condition $L\varphi\in\mathscr{G}$ is equal to the following 
\begin{equation}
  \label{eq:Lfi}
  L\varphi=f,
\end{equation}
where $f$ is uncertain and belongs to the given subset $\mathscr{G}$ of the
Hilbert space $\mathcal{F}$. In the sequel $\eta$ is supposed to be a random
$\mathcal{Y}$-valued vector with zero mean so that its correlation $R_\eta\in
\mathscr{R}$, where $\mathscr{R}$ is some subset
of $\mathscr{L}(\mathcal{Y},\mathcal{Y})$. Also we deal with deterministic
$\eta$ so that  $(f,\eta)\in\mathcal G$, where $\mathcal G$ is some subset of
$\mathcal F\times\mathcal Y$. Note that the realization of $y$ depends on
$\eta$, $H$ and $f$. Also it depends 
on elements of $N(L)=\{\varphi\in\mathscr D(L):L\varphi=0\}$ so that
$y=H(\varphi_0+\varphi)+\eta$, where $\varphi_0$ may be thought as inner noise
in the state model \eqref{eq:Lfi}. 

Let $\ell(\varphi)=(\ell,\varphi)$, $\widehat{\ell(\varphi)}=(u,y)+c$. Since
$L,H$ are not supposed to have a bounded inverse mappings the $\ell(\varphi)$
and $\widehat{\ell(\varphi)}$ are not stable with a respect to small
deviations in $f,\eta$. Also $f,\eta$ are supposed to be uncertain. Therefore
we use the minimax design in order to construct the estimation. 
\begin{dfn}\label{ozn1}
  The function $\widehat{\widehat{\ell(\varphi)}}=(\hat{u},\cdot)+\hat{c}$ is
  called the \emph{a priori minimax mean-squared estimation} iff
  $\sigma(\ell,\hat u)=\inf_{u,c}\sigma(\ell,u)$ where
\begin{equation}
    \label{eq:huc}
        \sigma(\ell,u):=\sup_{L\varphi\in \mathscr{G},R_\eta\in\mathscr R}
    M(\ell(\varphi)-\widehat{\ell(\varphi)})^2
  \end{equation}
The number $ \hat{\sigma}(\ell)=\sigma^\frac
12(\ell,\hat u)$ is said to be \emph{the minimax mean-squared error} in the
direction $\ell$. 
\end{dfn}
On the other hand the a posteriori estimation describes the evolution of the central point of the system reachability set $$
(L\varphi,y-H\varphi)\in\mathcal G
$$ consistent with measured output
$y$~\cite{Bertsekas1971,Tempo1985,Nakonechnii1978}. Note that the condition 
$(L\varphi,y-H\varphi)\in\mathcal G$ holds if $\|y\|<C$ for some real $C$. But
 it doesn't hold in our assumptions if $\eta$ is random since $\|R_\eta\|< c$
 doesn't imply $\|y\|<C$ for realizations of $\eta$. Therefore $\eta$ is
 supposed to be deterministic. 
\begin{dfn}\label{ozn2}
The set $$
\mathcal X_y=\{\varphi\in\mathscr D(L):(L\varphi,y-H\varphi)\in\mathcal G\}
$$ is called an a posteriori set. The vector $\hat{\varphi}$ is said to be
minimax a posteriori estimation of $\varphi$ in the direction $\ell$
($\ell$-minimax estimation) iff $$
\hat{d}(\ell):=\inf_{\varphi\in\mathcal X_y}\sup_{\psi\in\mathcal X_y}
|(\ell,\varphi)-(\ell,\psi)|=\sup_{\psi\in\mathcal X_y}
|(\ell,\hat{\varphi})-(\ell,\psi)|
$$ The expression $\hat{d}(\ell)$ is called the minimax a posteriori error in
the direction $\ell$ ($\ell$-minimax error).  
\end{dfn}
In the sequel the minimax mean-squared a priori estimation (error) is referred
as minimax estimation (error). 
\begin{prp}\label{t:1}
Assume that $\mathscr G$, $\mathscr{R}$ are convex bounded closed subsets of
$\mathcal F$, $\mathscr{L}(\mathcal{Y},\mathcal{Y})$ respectively. For the
given $\ell\in\mathcal H$ the minimax error $\hat\sigma(\ell)$ is finite iff 
\begin{equation}
  \label{eq:setUl}
  \ell-H^*u\in\mathrm{dom}\,\mathrm{cl}(L^*c)\cap
(-1)\mathrm{dom}\,\mathrm{cl}(L^*c)
\end{equation} for some $u\in\mathcal Y$. Under this condition 
  \begin{equation} 
    \label{eq:err:amap:thr}
    \begin{split}
          &\sigma(\ell,u)=\sup_{R_\eta\in\mathscr R}(R_\eta u,u)+\\
          &\frac 14[\mathrm{cl}(L^*c)(\ell-H^*u)+
    \mathrm{cl}(L^*c)(-\ell+H^*u)]^2
 \end{split}
    \end{equation}
where $$
R(L^*)\subset\mathrm{dom}\,\mathrm{cl}(L^*c)\subset\overline{R(L^*)}
$$
If $\mathrm{Arginf}_u\sigma(\ell,u)\ne\varnothing$, then
$\widehat{\widehat{\ell(\varphi)}}=(\hat{u},y)+\hat{c}$, where $$
\hat{u}\in\mathrm{Arginf}_u\sigma(\ell,u)$$
and $$
\hat{c}=\frac 12(\mathrm{cl}(L^*c)(\ell-H^*\hat u)-
\mathrm{cl}(L^*c)(-\ell+H^*\hat u))
$$
\end{prp}
\begin{thm}\label{t:RLintG}
  Suppose that $\mathscr G$ is convex bounded closed balanced set and $
  0\in\mathrm{int}\,\mathscr G$. Also assume that 
  $$
  \eta\in\{\eta:M(\eta,\eta)\le 1\}
  $$ Then for the given $\ell\in\mathcal H$ the minimax estimation
  $\hat\sigma(\ell)$ is finite iff $\ell-H^*u\in R(L^*)$  
  for some $u\in\mathcal Y$. Under this condition there exists a unique
  minimax estimation $\hat u$ and %\in\mathcal U_\ell$ and 
\begin{equation}
  \label{eq:aprer:amap:thr}
  \begin{split}
    &\sigma(\ell,\hat u)=\min_u\sigma(\ell,u),\\
    &\sigma(\ell,u)=(u,u)+\min_z\{c^2(\mathscr G,z),L^*z=\ell-H^*u\}
  \end{split}
\end{equation}
  If $R(L),H(N(L))$ are closed sets then $\hat u$ is determined by the
  following conditions
  \begin{equation}
    \label{eq:Hp0}
    \begin{split}
    &\hat u-Hp_0\in H(\partial I_2(H^*\hat u)),Lp_0=0,\\
    &I_2(w)=\min_z\{c^2(\mathscr G,z),L^*z=P_{L^*}(\ell-w)\},
    \end{split}
  \end{equation}
\end{thm}
\begin{cor}\label{t:2}
Let $$
\mathscr G=\{f\in\mathcal F:(f,f)\le 1\},
\eta\in\{\eta:M(\eta,\eta)\le 1\},
$$ and suppose that 
\begin{itemize}
\item [1)] $R(L),H(N(L))$ are closed sets;
\item [2)] $R(T)=\{[Lx,Hx],x\in\mathscr D(L)\}$ is closed set.
\end{itemize}
Then the unique minimax estimation $\hat u$
is given by $\hat u=H\hat p$ provided that $\ell\in R(L^*)+R(H^*)$, $\hat p$
obeys 
\begin{equation}
  \label{eq:qeuler}
  \begin{split}
    &L^*\hat z=\ell-H^*H\hat p,\\
    &L\hat p=\hat z
  \end{split}
\end{equation}
The minimax error is given by the following expression
\begin{equation*}
  %\label{eq:apr-err}
  \hat\sigma(\ell)=(\ell,\hat p)^\frac 12
\end{equation*}
\end{cor}
\begin{cor}\label{n:rozveuler}
Assume that linear mappings $
  L:\mathcal H\mapsto\mathcal F$, $
  H\in\mathscr L(\mathcal H,\mathcal Y)$ obey 1) or 2) (Cor.~\ref{t:2}).  
Then~\eqref{eq:qeuler} has a solution $
\hat z\in\mathscr D(L^*),\hat p\in\mathscr D(L)$ 
 iff $\ell=L^*z+H^*u$ for some $
z\in\mathscr D(L^*),u\in\mathcal Y$.
\end{cor}
\begin{cor}\label{n:alt}
Under conditions of Cor.~\ref{t:2} for any $\ell\in R(L^*)+R(H^*)$ and
some realization of $y(\cdot)$ we have $(\hat u,y)=
(\ell,\hat\varphi)$, where $\hat\varphi$ obeys 
\begin{equation}
  \label{eq:alt}
  \begin{split}
    &L^*\hat q=H^*(y-H\hat\varphi),\\
    &L\hat\varphi=\hat q
  \end{split}
\end{equation}
\end{cor}
Consider an a posteriori estimation. 
\begin{prp}\label{t:apo:1}
Let $\mathcal G$ be a convex closed bounded subset of $\mathcal
  Y\times\mathcal F$. Then 
  \begin{equation}
    \label{eq:setL}
    \begin{split}
      &R(L^*)+R(H^*)\subset\mathrm{dom}\,c(\mathcal X_y,\cdot)\cap
    (-1)\mathrm{dom}\,c(\mathcal X_y,\cdot)\subset\\
    &\overline{R(L^*)+R(H^*)}
    \end{split}
  \end{equation}
The minimax a posteriori error in the direction $\ell$
  is finite iff $\ell\in
  \mathrm{dom}\,c(\mathcal X_y,\cdot)\cap(-1) 
  \mathrm{dom}\,c(\mathcal X_y,\cdot)$ and
  \begin{equation}
    \label{eq:apo:esterr}
    \begin{split}
      &(\ell,\hat\varphi)=\frac 12(c(\mathcal X_y,\ell)-c(\mathcal
    X_y,-\ell)),\\
    &\hat{d}(\ell)=\frac 12(c(\mathcal X_y,\ell)+c(\mathcal X_y,-\ell))
    \end{split}
  \end{equation}
\end{prp}
\begin{thm}\label{t:apo:2}
Let $$
\mathcal G=\{(f,\eta):\|f\|^2+\|\eta\|^2\le 1\},
$$ and assume that $1)$ or $2)$ from Corollary~\ref{t:2} holds. The
minimax a posteriori estimation $\hat{\varphi}$ obeys %the equation 
\begin{equation}
  \label{eq:apoest:amap:thr}
  \begin{split}
    &L^*\hat q=H^*(y-H\hat\varphi),\\
    &L\hat\varphi=\hat q
  \end{split}
\end{equation}
iff $\ell\in R(L^*)+R(H^*)$. The estimation error is given by 
\begin{equation}
  \label{eq:apos-err}
  \hat{d}(\ell)=(1-(y,y-H\hat\varphi))^\frac 12 \hat\sigma(\ell)
  %(\ell,\hat p)^\frac 12
\end{equation}
\end{thm}
\begin{cor}\label{amap:thr:vest}
  Assume that the conditions of Theorem~\ref{t:apo:2} are fulfilled and
  $\widehat{\ell(\varphi)}=(\ell,\hat\varphi)$ for any $\ell$, where
  $\hat{\varphi}$ obeys~\eqref{eq:apoest:amap:thr}. Then 
  %$\hat\varphi$ gives
  %the minimax a posteriori estimation of $\varphi$ so that 
\begin{equation*}
  \begin{split}
    &\inf_{\varphi\in\mathcal X_y}\sup_{x\in\mathcal X_y}\|\varphi-x\|=\\
    &\sup_{x\in\mathcal X_y}\|\hat\varphi-x\|=(1-(y,y-H\hat\varphi)^\frac 12
    \max_{\|\ell\|=1}\hat{\sigma}(\ell)
    %\|(L^*L+H^*H)^{-\frac 12}\|\times\\
    %&%(1-((-H(L^*L+H^*H)^{-1}H^*)y,y)^\frac 12
  \end{split}
\end{equation*}
\end{cor}
In order to apply these results for linear DAEs we investigate some properties
of the linear mapping induced by DAE \cite{Zhuk2007}. 
Let $$
L\varphi(t)=[\dfrac d{dt}
F\varphi(t)-C(t)\varphi(t),F\varphi(t_0)]
$$ and set $$
\mathscr
D(L)=W_F:=\{\varphi(\cdot)\in \mathbb L_2(t_0,T):t\mapsto F\varphi(t)\in \mathbb
W_2(t_0,T)\}
$$ 
It is clear that $L\varphi(t)=[f(t),f_0]$ is equal to $$
\dfrac d{dt}F\varphi(t)=C(t)\varphi(t)+f(t), F\varphi(t_0)=f_0
$$ Next proposition describes the adjoint $L^*$. 
\begin{thm}%\label{t:3}
If $x(\cdot)\in W_F$, $z\in W_{F'}$ then \begin{equation}
  \label{eq:ich}
  \begin{split}
  &\int_{t_0}^T(\dfrac d{dt}Fx(t),z(t))+(\dfrac d{dt} F'z(t),x(t))dt=\\
  &(Fx(T),F'^+F'z(T))-(Fx(t_0),F'^+F'z(t_0))
\end{split}
\end{equation}
$L$ is closed linear mapping and its adjoint $L^*:\mathbb L_2(t_0,T)\times
\mathbb R^n\to \mathbb L_2(t_0,T)$ is defined as follows 
\begin{equation*}
  %\label{eq:Ds}
  \begin{split}
  &L^*(z,z_0)(t)=-\dfrac d{dt} F'z(t)-C'(t)z(t),\\
  &\mathscr{D}(L^*)=\{(z,F'^+F'z(t_0)+d):z\in W_{F'},F'z(T)=0,F'd=0\}
\end{split}
\end{equation*}
\end{thm}
Note that $R(L)$ is not necessary close. A sufficient condition
for $R(L)$ to be close is introduced in the next theorem
assuming\footnote{This assumption holds for any linear DAE with constant
  matrices. } that $$
F=\bigl(
\begin{smallmatrix}
  E&&0\\0&&0
\end{smallmatrix}
\bigr),
C=\bigl(
\begin{smallmatrix}
  C_1&&C_2\\C_3&&C_4
\end{smallmatrix}
\bigr)
$$
\begin{thm}
If $$
\sup_{1>\varepsilon>-1}\|Q(\varepsilon)C'_2\|_{mod}<+\infty,Q(\varepsilon):=(\varepsilon^2
E+C'_4C_4)^{-1}, \|F\|_{mod}:=\sum_{i,j}|F_{ij}|
$$ then $R(L)$ is closed. 
\end{thm}

Next subsections demonstrate the application of the above theory to the linear
estimation problem for linear DAEs.   
\subsection{DAEs with continuous time}
In this subsection we present the main result of preprint \cite{Zhuk2009b}
-- linear reduced order minimax filter for linear noncausal DAEs with continuous
time. All proofs are given in \cite{Zhuk2009b}. 

Consider a pair of systems% of linear differential-algebraic equations
  \begin{equation}
    \label{eq:dae}
    \begin{split}
      & \dfrac d{dt}Fx(t)=C(t)x(t)+f(t),Fx(t_0)=0,\\
      & y(t)=H(t)x(t)+\eta(t),t\in[t_0,T],
    \end{split}
   \end{equation}
where $x(t)\in\mathbb R^n$, $f(t)\in\mathbb R^m$, $y(t)\in\mathbb R^p$,
$\eta(t)\in\mathbb R^p$ represent the state, input, measurement output and
measurement noise respectively, $F\in\mathbb R^{m\times n}$, $t\mapsto
C(t)\in\mathrm C^{m\times n}(t_0,T)$, $f(\cdot)\in\mathbb L_2(t_0,T)$,$t\mapsto H(t)\in\mathrm C^{p\times n}(t_0,T)$, $t_0,T\in\mathbb R$.  

According to~\cite{Zhuk2007} we say that $x(\cdot)$ is a
solution of~\eqref{eq:dae} if $Fx(\cdot)\in\mathbb W_2^m(t_0,T)$, the
derivative of $Fx(\cdot)$ coincides with the right side
of~\eqref{eq:dae} almost everywhere and $Fx(t_0)=0$ holds. 

In the sequel we assume that $\eta(\cdot)$ is
a realization of the random process $\eta$ with zero mean satisfying  
\begin{equation}
  \label{eq:eta_bounds}
\eta\in W=\{\eta:M\int_{t_0}^T(R(t)\eta(t),\eta(t))\le 1  \}
\end{equation}
and $$
f(\cdot)\in G=\{f(\cdot):\int_{t_0}^T(Q(t)f(t),f(t))\le 1\},
$$
where $Q(t)\in\mathbb R^{m\times m}$, $Q=Q'>0$, $R(t)\in\mathbb
R^{p\times p}$, $R'=R>0$ and $Q(t),R(t)$, $R^{-1}(t)$, $Q^{-1}(t)$ are
continuous functions of $t$ on $[t_0,T]$. 

Suppose $y(t)$ is observed in \eqref{eq:dae} for some $x(\cdot)$, $f\in G$ and
$\eta$. Our aim here is to construct an algorithm\footnote{In literature it is
common to refer to this algorithm as filter~\cite{Brammer1989}} giving
online estimation of the linear function $$
x(\cdot)\mapsto (\ell, Fx(T))
$$ on the basis of the measured on $[t_0,T]$ realization of the output
$y(t)$. 
With this purpose we
introduce a notion of the linear minimax estimation~\cite{Zhuk2008b}.    
\begin{defn}
  The function $\hat u(y)=\int_{t_0}^{T}(\hat u(t),y(t))dt$ is called
minimax mean-squared a priori estimation if $$ 
\inf_u\sigma(u)=\sigma(\hat u)
$$ where $
\sigma(u)=\sup_{x(\cdot),f(\cdot),\eta}M[(\ell, Fx(T))-u(y)]^2
$ is a maximum estimation error for $u(\cdot)$. The number
$\hat\sigma=\sigma(\hat u)$ is called a minimax mean-squared a priori
error. The state $x(t)$ is said to be minimax observable in the direction
$\ell$ iff $\hat\sigma<+\infty$. 
\end{defn}
The minimax
directional observability differs from the classical observability property in
the following way. If system state $x(s)$ is minimax observable in the
direction $\ell$ then the projection of the reachability set (consistent with
some realization of $y(t),t_0\le t\le s$) onto direction $\ell$ is expected to
be $[-\hat\sigma,\hat\sigma]$, where $\hat\sigma$ denotes the minimax estimation
error. The expected estimation error varies in $[0,\hat\sigma]$ and depends on
the noise realization, initial condition and input. If $x(s)$ is unobservable in
the minimax sense for $\ell$ then the minimax estimation is set to zero and
$\hat\sigma=+\infty$. This means that the structure of the measurements do not
provide any information about $(\ell,x(s))$. Therefore the minimax
directional observability provides a qualitative description of DAEs
singularity with a respect to the given observations. In particular, regular
DAE is observable in the minimax sense for any direction in contrast to
the classical observability.
% If we consider the set of numbers $\{r\in[0,+\infty]:r=g(x,y)\}$, where $(x,y)$
% runs through the set of all possible functions that satisfy~\eqref{eq:dae}
% for any $f\in G$ and any realization of $\eta\in W$, then the symbol
% $\sup_{x(\cdot),f(\cdot),\eta}g(x,y)$ denotes its upper bound.
% In the sequel the minimax mean-squared a priori estimation (error) is referred
% as minimax estimation (error), the argument $t$ is dropped almost everywhere.
\begin{rem}
  This definition generalizes the notion of linear minimax a priori estimation
  introduced in \cite{Nakonechnii1978}. Here we follow a common way of  
  \cite{Nakonechnii1978} deriving the minimax 
  estimation: first step is to describe a dual control problem, next step is
  to solve it and the last step is to derive a minimax filter.  
\end{rem}

Assume that $u(\cdot)\in\mathbb L_2(t_0,T)$, $\ell\in\mathbb R^m$ and $z(\cdot)$
obeys DAE %denotes any solution of the boundary value problem (BVP)
\begin{equation}
    \label{eq:zul}
    \dfrac d{dt}F'z(t)=-C'(t)z(t)+H'(t)u(t),F'z(T)=F'\ell
  \end{equation}
Let $v(\cdot)$ denotes any solution of homogeneous DAE \eqref{eq:zul}. Next
proposition gives a generalization of the celebrated Kalman duality
principle~\cite{Brammer1989}.    
\begin{prop}\label{p:1}
The minimax estimation error $$
\sigma(u)=\sup_{x(\cdot),f(\cdot),\eta}M[(\ell, Fx(T))-u(y)]^2\to\inf_u
$$ is finite iff \eqref{eq:zul} has a solution $z(\cdot)$. The minimax
estimation problem $\sigma(u)\to\inf_u$ is equal to the following optimal
control problem 
\begin{equation}
  \label{eq:umin}
  I(u)=\min_v\{\int_{t_0}^T(Q^{-1}(z-v),z-v)dt\}+\int_{t_0}^T(R^{-1}u,u)dt\to\min_u,
\end{equation}
provided that $z(\cdot)$ is some solution of \eqref{eq:zul}. 
% where \begin{equation}
%   \label{eq:NL*}
% \dfrac d{dt}F'v(t)=-C'(t)v(t),F'v(T)=0
% \end{equation}
\end{prop}
Proposition~\ref{p:1} states that minimax estimation problem is equal to some
optimal control problem for appropriate $\ell$ which is called dual control
problem. In the next proposition we
introduce a representation for the minimax estimation and error. 
\begin{prop}\label{p:2}
  Let $p(\cdot)$ denotes some solution of the two-point boundary value problem
  \begin{equation}
    \label{eq:dae_bvp}
    \begin{split}
      &\dfrac d{dt}Fx(t)=C(t)x(t)+Q^{-1}(t)z(t),Fx(t_0)=0,\\
      &\dfrac d{dt}F'z(t)=-C'(t)z(t)+H'(t)R(t)H(t)p(t),F'z(T)=F'\ell
    \end{split}
  \end{equation}
Then minimax estimation $\hat u$ is given by $\hat u=RHp$, the minimax error
is represented as $\hat\sigma=(\ell,Fp(T))$. 
\end{prop}
It is known that ~\eqref{eq:dae} may be converted into SVD coordinate system
\cite{Bender1987} so that without loss of generality we assume that $$
F=\bigl(
\begin{smallmatrix}
  E&&0\\0&&0
\end{smallmatrix}
\bigr),
C=\bigl(
\begin{smallmatrix}
  C_1&&C_2\\C_3&&C_4
\end{smallmatrix}
\bigr),Q(t)=\bigl(
\begin{smallmatrix}
  Q_1&&Q_2\\Q_3&&Q_4
\end{smallmatrix}
\bigr),
R(t)=\bigl(
\begin{smallmatrix}
  R_1&&R_2\\R_3&&R_4
\end{smallmatrix}
\bigr),S(t)=\bigl(
\begin{smallmatrix}
  S_1&&S_2\\S_3&&S_4
\end{smallmatrix}
\bigr)
$$ where $S=H'RH$. By definition, put $A(t)=C_1-Q_2Q^{-1}_4C_3-(C_2-Q_2Q^{-1}_4C_4)\tilde
S^+_4(S_3+C_4'Q^{-1}_4C_3)$,
$M(t)=S_1+C_3'Q_4^{-1}C_3-(S_2C_3'Q_4^{-1}C_4)\tilde
S^{+}_4(S_3+C_4'Q_4^{-1}C_3)$, $\overline C=\tilde
S^{+}_4((C_2-Q_2Q^{-1}_4C_4)^T-(S_3+C_4'Q^{-1}_4C_3)K)$,
$G(t)=Q_1-Q_2Q_1^{-1}Q_3+(C_2-Q_2Q_4^{-1}C_4)\tilde
S_4^+(C_2-Q_2Q_4^{-1}C_4)^T$.
\begin{thm}
  Assume that $t\mapsto \tilde S_4^+(t)=(S_4(t)+C_4'(t)Q_4^{-1}(t)C_4(t))$ is
  measurable matrix-valued function. For any $\ell\in\mathbb R^n$ the minimax
  estimation of the inner product $(\ell,Fx(T))$ is given by $$ 
\widehat{(\ell,Fx(T))}=(\ell_1,\hat x(T))
$$ where $\hat x$ is the solution of the initial-value problem
\begin{equation}
  \begin{split}
    &\dfrac d{dt}\hat x=(A(t)-K(t)M(t))\hat x+K(t)[E,\overline C']H'Ry(t), 
    \hat x(t_0)=0,\\
    &\dot K=AK+KA'+KMK-G,K(t_0)=0
  \end{split}
\end{equation}
The minimax estimation error is given by $\hat \sigma =(\ell_1,K(T)\ell_1)$,
where $\ell$ is splitted into $(\ell_1,\ell_2)$ according to the block
structure of $F$.  
\end{thm}
\subsection{DAEs with discrete time}
In this subsection we present the main result of the preprint \cite{Zhuk2009a}
-- linear recursive minimax filter for linear noncausal DAE with discrete
time.
 All proofs are given in \cite{Zhuk2009a}.  

Consider %a LTV system described by 
the %following discrete-time descriptor 
model
\begin{align}
  %  \begin{split}
    &F_{k+1}x_{k+1}-C_k x_k = f_k, F_0x_0=q,\label{eq:state}\\
    &y_k=H_k x_k+g_k, k=0,1,\dots\label{eq:mes}
%  \end{split}
\end{align}
where $x_k\in\mathbb R^n$, $f_k\in\mathbb R^m$, $y_k,g_k\in \mathbb R^p$ % and
% $g_k\in\mathbb R^p$
represent the state, input, measurement output and
measurement noise respectively, $F_k,C_k\in\mathbb R^{m\times n}$, $H_k\in
\mathbb R^{p\times n}$ % for $k=0,1,\dots,$
and initial state $x_0$ belongs to
the affine set $\{x:F_0x=q\}$, $q\in\mathbb R^m$.
%, $F_0\in\mathbb R^{s\times n}$.  
% \begin{equation}
%   \label{eq:strt}
% \end{equation}
% and $y_k$ is given by 
% \begin{equation}
%   \label{eq:obs}
%   % k=\overline{0,\uN},
% \end{equation}
% where , $H_k$ is $p\times n$-matrix. 
In what follows we assume that  
%[q,\{f_s\}_0^\infty,]$ 
%We assume that 
\begin{equation}
  \label{eq:ellips}
\xi\in\mathscr G=
\{\xi:\Psi_n(\xi_n)%=(S q,q)+ \sum_0^{n-1}(S_i f_i,f_i)+(R_ig_i,g_i)
\leqslant1,\forall n\in\mathbb N\}  
\end{equation}
where $\xi=[q,\{f_s\},\{g_s\}]$, $\{f_s\}=[f_0,f_2,\dots]$, 
$\xi_k=[q,\{f_s\}_0^k,\{g_s\}_0^k]$, $\{f_s\}_0^k=[f_0\dots f_k]$ is the
projection of $\{f_s\}$ onto linear span of $e_1\dots e_{k+1}$,
$e_1=[1,0\dots]$, $ 
  \Psi_n(\xi_n)=(S q,q)+ \sum_0^{n-1}(S_i f_i,f_i)+(R_i g_i,g_i)
$, $S$, $S_k\in\mathbb R^{m\times m}$, $R_k\in \mathbb R^{p\times p}$ are
symmetric and positive-definite. 

%, $S=S'>0$, $S_k=S'_k>0$, $R_k=R'_k>0$.\\
 Suppose that $y^*_k$ %$k=0,1,\dots$ 
 is being observed in~\eqref{eq:mes} with
 $x_k=x^*_k$ and $g_k=g^*_k$ provided that $x^*_k$ is derived from
 \eqref{eq:state} with $f_k=f^*_k$, $q=q^*$ and
 $\xi^*=[q^*,\{f^*_s\},\{g^*_s\}]\in\mathscr G$. 
\emph{Our aim here is} to describe the evolution in $\tau$ of the $\ell$-minimax
estimation $\overline x_\tau$ of the state $x^*_\tau$ at instant $k=\tau$ along
with error $\hat\rho(\ell,\tau)$ through dynamic
recurrence-type relation and to efficiently describe the structure of the
minimax observable subspace $\mathcal L(\tau)$. For this purpose we shall
apply the theory developed in the previous section.
\begin{defn}\label{d:2}
The set $\mathcal L(\tau)=\{\ell:\hat{\rho}(\ell,\tau)<\infty\}$ is
called a minimax observable subspace for the model~\eqref{eq:state} at the
instant $k=\tau$. Its co-dimension $I_\tau=n-\mathrm{rank Q_\tau}$ is called
\emph{an index of non-causality} of the model~\eqref{eq:state}.
\end{defn}
\begin{thm}\label{t:mnmx}%[minimax recursive estimation]
  The minimax observable subspace for the model~\eqref{eq:state} is given
  by $\mathcal L(\tau)=\{\ell:P_\tau^+P_\tau\ell=\ell\}$ and
  \begin{equation}
    \label{eq:Xtau}
    X(\tau)=\{x\in\mathbb R^n:(P_\tau(x-\hat x_\tau),x-\hat x_\tau)\le\hat\beta_\tau\}
 \end{equation}
 where  $\hat x_\tau=P_\tau^+r_\tau$, $\hat\beta_\tau:=1-\alpha_\tau+(P_\tau\hat x_\tau,\hat x_\tau)$ 
\begin{equation*}
    %\label{eq:Qk}
    \begin{split}
      &P_k=%F'_kS_{k-1}F_k+
      H'_kR_kH_k+F'_k[S_{k-1}-S_{k-1}C_{k-1}B_{k-1}^+C'_{k-1}S_{k-1}]
      F_k,\\
      &P_0=F'_0SF_0+H'_0R_0H_0,B_k=P_{k}+C'_{k}S_{k}C_{k},
    \end{split}
  \end{equation*}  
$\alpha_i=\alpha_{i-1}+(R_iy_i,y_i)-(B_{i-1}^+
 r_{i-1},r_{i-1})$, $\alpha_0=(R_0y_0,y_0)$ and $$
r_k=F'_kS_{k-1}C_{k-1}B_{k-1}^+r_{k-1}+H'_kR_k y_k,
      r_0=H'_0R_0y_0
$$ 
If $\ell\in\mathcal L(\tau)$ then $\hat x_\tau$ is
$\ell$-minimax estimation of $x^*_\tau$ and $
    \hat{\rho}(\ell,\tau)={\hat\beta}_\tau^\frac12
    %[1-\alpha_\tau+(P^+_\tau r_\tau,r_\tau)]^\frac 12
    (P^+_\tau\ell,\ell)^\frac 12
$.  
%If $I_\tau=0$ then    
\end{thm}
\textbf{Example.} %\label{s:exmpl}
% \textbf{ToDo:stress the benefits of non-causal state estimation!$v_k(x_k)$ --
%   feedback control! }
 Consider the following filtration problem: given measurements $y_k$ one needs
 to construct the estimation $\hat x_\tau$ of the state $x_k$ at instant $k=\tau$ and to describe the estimation error provided that
 %$$
 \begin{equation}
   \label{eq:exm1}
   x_{k+1}=c_kx_k+v_k(x_k)+f_k,x_0=q,  
y_k=h_kx_k+w_k%\eqno(*)
 \end{equation}
%$$ 
%$v_k$, $q$, $f_k$ and $w_k$ 
$p\mapsto v_k(p)$ is some real-valued function and  
uncertain scalar parameters $q$, $w_k$ and $f_k$ are
restricted by the inequality 
\begin{equation}
  \label{eq:exG}
  %$$%I(x_0,w_k,f_k)=
Sq^2+\sum_{k=0}^{\tau-1}R_kw_k^2+S_kf^2_k\le 1, S,R_k,S_k>0%\eqno(**)
\end{equation}%$$
%We stress that an a priori knowledge about $v_k$ is not available. In
% practice it may happen when $v_k$ corresponds to an exogenious unmeasured
% control parameter.
Let us show how one can construct $\hat x_\tau$ 
%so that the worst-case error
%do not depend on $v_k$ realization. This will be done 
by means of the set-membership state estimation approach for linear non-causal
descriptor systems described in this section. %Theorem 
Let $z_k=[z_{1,k},z_{2,k}]$ obeys %$$
\begin{equation}
  \label{eq:exFz}
  Fz_{k+1}=C_kz_k+f_k,\, Fz_0=q,\, 
a_k=H_kz_k+w_k %\eqno(***)
\end{equation}%$$ 
where $F=(1,0)$, $C_k=(c_k,1)$,
$H_k=(h_k,0)$. Note that for any real
$z_{2,k}$ there exists exactly one $z_{1,k}$ %(derived from ~\eqref{eq:exFz})
so that $z_k$ obeys the first equation in~\eqref{eq:exFz}. \\
We shall apply Theorem~\ref{t:mnmx} in order to construct the
$\ell$-minimax estimation of $z_\tau$. Using definitions of $P_k,r_k$ one obtains $P_0=
 \left(\begin{smallmatrix}
   q_0&&0\\0&&0
 \end{smallmatrix}\right)
$, $r_0=R_0H'_0a_0$, where $q_0=S+R_0$ so that $B^+_0=%(Q_0+C'_0S_0C_0)^+=
\left(
  \begin{smallmatrix}
    \frac{1}{q_0}&&-\frac{c_0}{q_0}\\
    -\frac{c_0}{q_0}&&\frac{c_0^2}{q_0}+\frac{1}{S_0}
  \end{smallmatrix}\right)
$ and therefore $S_0C_0B^+_0=(0,1)$, $S_0C_0B^+_0C'_0=1$ so that
$P_1=R_1H'_1H_1$ and $r_1=R_1H'_1a_1$. It's easy to prove by induction that
$P_k=R_kH'_kH_k$ and $r_k=R_kH'_ka_k$. Theorem~\ref{t:mnmx} implies: $\hat
z_\tau:=P^+_\tau r_\tau$ represents the $\ell$-minimax estimation for
$\ell\in\mathcal L(\tau)$, where  
$\mathcal L(\tau)=\{\lambda e_1,e_1=[1,0],\lambda\in\mathbb R\}$ if
$h_\tau\ne0$ and $\mathcal L(\tau)=\{0\}$ otherwise. Hence if $\ell=[0,l]$
then %$\ell\not\in\mathcal L(\tau)$ so that 
the $\ell$-minimax
error is infinite. This fact reflects the non-causality
of the model \eqref{eq:exFz}: since the a
 posteriori set of \eqref{eq:exFz} is a shift of
 the convex set $
\mathscr G(0)=\{[z_0\dots z_\tau]:$
$$\|
\bigl(\begin{smallmatrix}
  F\\H
\end{smallmatrix}\bigr)z_0\|^2+
%Fz_0\|^2+\|Hz_0\|^2+
\sum_0^{\tau-1}\|\bigl(\begin{smallmatrix}
  F&&-C_k\\H_{k+1}&&0
\end{smallmatrix}\bigr)
\bigr[\begin{smallmatrix}
  z_{k+1}\\z_k
\end{smallmatrix}\bigl]\|^2\le\beta,\beta>0\}
%Fz_{k+1}-Bz_{k-1}\|^2+\|Dz_k\|^2\le 10\}
$$ its cross section $Z(\tau)$ at the instant $k=\tau$ is a shift of $\mathcal
P_\tau(\mathscr G(0))$. Thus $Z(\tau)$ is convex and unbounded implying that
it recedes to infinity \cite[\S 8]{Rockafellar1970} in the directions
$\ell\notin\mathcal L(\tau)$.  
% This fact reflects the non-causality of~\eqref{eq:exFz}:
% since $z^2_k$ is not restricted to belong to some bounded set the cross
% section of the a posteriori set of~\eqref{eq:exFz} is unbounded in the direction
% $\ell=[0,l]$. 
%But its projection onto minimax observable subspace $\mathcal L$
%is bounded. This is reflected in relation~\eqref{eq:exz1t}.Therefore $I_\tau\le 1$ so
%that~\eqref{eq:exFz} is non-causal. 
If $h_\tau\ne 0$ and $\ell=[l,0]$ then the 
$\ell$-minimax estimation $\hat z_\tau$ %and error $\hat\rho^2(\ell,\tau)$ 
obeys
\begin{equation}
  \label{eq:exmErEst}
  (\ell,\hat z_\tau)=\frac{l}{h_\tau}a_\tau=l(z^1_\tau+\frac{w_\tau}{h_\tau}),
(\ell,\hat z_\tau-z_\tau)^2%\hat\rho^2(\ell,\tau)
\le\frac{l^2}{R_\tau h^2_\tau}
\end{equation}
since $w_\tau^2\le R_\tau^{-1}$ due to~\eqref{eq:exm1}. 
%Note that \eqref{eq:exmErEst} implies % $\ell$-minimax error $\hat\rho^2(\ell,\tau)$ doesn't depend on $S_\tau$. This in turn means that one can
%  choose $S_k$ sufficiently small in order to take into account realizations of
%  $\{f_k\}_{k\le\tau}$ with large energy. Also \eqref{eq:exmErEst} implies
% that the $\ell$-minimax error decays iff $h^2_kR_k$ grows. 
% If $\ell\in\mathcal
%  L(\tau)$ then 
% $$(\ell,z^\tau)=\frac{\ell_1}{h_\tau}y^1_\tau=\ell_1(z^1_\tau+\frac{w_\tau}{h_\tau})$$
%\newline 
Let $y^*_k,x^*_k$ denote
the realization of output $y_k$ and state $x_k$ derived from \eqref{eq:exm1}
 with $f_k=f^*_k$, $q=q^*$, $w_k=w^*_k$ % and $v_k=v^*_k$. 
%We assume that \eqref{eq:exG} holds and $v^*_k$ is some real-valued function of the
%natural parameter $k$. 
restricted by \eqref{eq:exG}. Let $z^*_{k}=[z^{*}_{1,k},z^{*}_{2,k}]$, $a^*_k$
be derived from~\eqref{eq:exFz} provided that $f_k=f^*_k$,
$z^{*}_{2,k}:=v_k(z^*_{1,k})$, $q=q^*$ and $w_k=w^*_k$. By direct calculation $z^*_{1,k}=x^*_k$ and $a^*_k=y_k^*$ so that
$lx^*_\tau=(\ell,z^*_\tau)$ with $\ell=[l,0]$. Hence 
%and $z_\tau^*=[z^{1,*}_\tau,z^{2,*}_\tau]$. Hence
\begin{equation*}
%  \label{eq:exz1t}
  (lx_\tau^*-(\ell,\hat z_\tau))^2=\bigl(\frac{lw^*_\tau}{h_\tau}\bigr)^2\le
\frac{l^2}{R_\tau   h^2_\tau}%=\hat\rho^2(\ell,\tau)
\end{equation*} 
due to~\eqref{eq:exmErEst}. Thus %the function
$\tau\mapsto(\ell,\hat z_\tau)$ gives the online estimation of $\tau\mapsto
lx_\tau^*$ with worst-case error $\frac{l^2}{R_\tau   h^2_\tau}$.
%described by \eqref{eq:exz1t}.\\
\newline\textbf{Minimax estimator and $H_2/H_\infty$ filters.}
In \cite{Baras1995} a connection between set-membership state estimation
and $H_\infty$ approach is described for linear causal DAEs. The authors note
that the notion of informational state ($X(\tau)$ in our notation) is shown to
be intrinsic for both approaches: the mathematical relations between
informational states of $H_\infty$ and set-membership state estimation are
described in \cite[Lemma 6.2.]{Baras1995}. Comparisons of set-membership
estimators with $H_\infty$ and other widely used filters for linear
DAEs are presented in \cite{Sayed2001} provided that $F_k\equiv E$. \\ 
In~\cite{ishixara2} % Kalman's filtering problem for descriptor systems was
%   investigated from the deterministic point of view.
  authors recover Kalman's
  recursion to LTV DAE from a deterministic least
  square fitting problem over the entire trajectory: %:  
%\begin{flushleft}
%(\textrm{DLSFP}) 
% find a sequence $\{\hat{x}_{0|k},\dots,\hat{x}_{k|k}\}$ that minimises the following fitting error cost
% \begin{equation*}
%   \begin{split}
%     &\mathrm{J}_k(\{x_{i|k}\}_{0}^k)=
% \|F_0x_{0|k}-g\|^2+\|y_0-H_0x_{0|k}\|^2+\\
% &\sum_{i=1}^k\|F_{i}x_{i|k}-C_{i-1} x_{i-1|k}\|^2+\|y_i-H_ix_{i|k}\|^2
%   \end{split}
% \end{equation*}
% assuming that $\mathop{\mathrm{rank}}
%   \begin{smallmatrix}
%     F_k\\H_k
%   \end{smallmatrix}\equiv n$. 
%\end{flushleft}
%According to~\cite{ishixara2} 
if $\mathop{\mathrm{rank}}
  \begin{smallmatrix}
    F_k\\H_k
  \end{smallmatrix}\equiv n$ then the optimal estimation $\hat{x}_{i|k}$ % resulting from the minimisation of $\mathrm{J}_k$
can be found from %the recursive algorithm 
  \begin{equation*}
    %\label{eq:fltr:r}
    \begin{split}
     &\hat{x}_{k|k}=
   P_{k|k}F'_{k}A_{k-1}C_{k-1}\hat{x}_{k-1|k-1}+P_{k|k}H'_{k} R_{k}y_{k},\\
   &\hat{x}_{0|0}=P_{0|0}H'_0R_0y_0,%P_{0|0}(F'_0Sq+H'_0R_0y_0),
   A^{-1}_k=(S^{-1}_k+C_{k}P_{k|k}C'_{k})\\
  &P^{-1}_{k|k}=F'_kA_{k-1}F_k+H'_kR_kH_k,
  P^{-1}_{0|0}=F'_0SF_0+H'_0R_0H_0
    \end{split}
  \end{equation*}
%\end{rem}
\begin{cor}
\label{l:eqfltr}  Let $r_0=F'_0Sq+H'_0R_0y_0$. If $\mathop{\mathrm{rank}}
  \bigl[\begin{smallmatrix}
    F_k\\H_k
  \end{smallmatrix}\bigr]\equiv n$  
%, and update $ r_k$  according to~\eqref{eq:rk1}. %  be a recursive map that
%   takes each natural number $k$ to the vector $r_k\in\mathbb R^n$, where 
% \begin{equation}
%     \label{eq:rk2}
%     \begin{split}
%       &r_k=H'_k y_k+F'_k C_{k-1}(C'_{k-1}C_{k-1}+Q_{k-1})^+_{k-1}r_{k-1},\\
%       &r_0=F'_0q+H'_0y_0
%     \end{split}
% \end{equation}
then $I_k=0$ and $P^+_kr_k=\hat{x}_{k|k}$.% for each $k\in\mathbb{N}$.% , where
% $\hat{x}_{k|k}$ is given by~\eqref{eq:fltr:r}.   
\end{cor}
\section{Numerical example}\label{s:exmpl}
Let us show how to use the minimax estimation in the
infinite-horizon setting. Consider the following DAE
\begin{align}
  %  \begin{split}
    &\left[
\begin{smallmatrix}
  1&&0&&0\\0&&k&&0
\end{smallmatrix}\right]\left(
\begin{smallmatrix}
x_{1,k+1}\\
x_{2,k+1}\\
x_{3,k+1}
\end{smallmatrix}\right)=
\left[
\begin{smallmatrix}
  \frac 1{40}&&\frac 12&&0\\\frac 1{10}&&\frac 14&&\frac 3{10}
\end{smallmatrix}\right]\left(
\begin{smallmatrix}
x_{1,k}\\
x_{2,k}\\
x_{3,k}
\end{smallmatrix}\right)+\left(
\begin{smallmatrix}
f_{1,k}\\
f_{2,k}
\end{smallmatrix}\right)\label{eq:st:ne}\\
%    &kx_{2,k+1}=\frac {x_{1,k}}{10}+\frac{x_{2,k}}4+\frac{x_{3,k}}{10}+f^k(r_2,p_2),\\
    &\left(
\begin{smallmatrix}
y_{1,k}\\
y_{2,k}\\
y_{3,k}\\
y_{4,k}
\end{smallmatrix}\right)=\left[ 
\begin{smallmatrix}
  h_{1,k}&&h_{2,k}&&0\\
  h_{4,k}&&h_{5,k}&&0\\
  h_{8,k}&&0.005&&h_{3,k}\\
  h_{6,k}&&h_{7,k}&&0
\end{smallmatrix}\right]\left(
\begin{smallmatrix}
x_{1,k}\\
x_{2,k}\\
x_{3,k}
\end{smallmatrix}\right)+
\left(
\begin{smallmatrix}
g_{1,k}\\
g_{2,k}\\
g_{3,k}\\
g_{4,k}
\end{smallmatrix}\right)\label{eq:obs:ne}
\end{align}
%  \end{split}
%\end{equation}
where % $f_{1,k}=f_1^k(p_1)$, $f_{2,k}=f_2^{k+300}(p_2)$,
% $g_{1,k}=f_1^{k+400}(p_1)$, $g_{2,k}=f_2^k(p_2)+f_1^k(p_1)$,
% $g_{3,k}=\sin\exp(k)$, $g_{4,k}=f_3^k(p_3)$,
$h_{1,k}=\frac {6k}{10}$,
$h_{1,0}=\frac {6}{10}$, $h_{3,k}=150k$ if $k$ is odd and $0$ otherwise;
$h_{4,k}=100k$, $h_{4,0}=1000$, $h_{2,k}=k$, $h_{2,0}= \frac {96}{100}$,
$h_{5,k}=\frac k{100}$, $h_{5,0}=2\frac 3{10}$, $ h_{6,k}=0.05$, $h_{6,0}=0$,
$h_{7,k}=10k$, $h_{7,0}=0$, $h_{8,k}=0$, $h_{8,0}=1$. Also we set
$x_{1,0}=1,x_{2,0}=-3$ and suppose that 
$$
(Sq,q)+\sum_0^\infty(R_kg_k,g_k)+(S_kf_k,f_k)\le 1
$$
where  $
R_k=\frac 1{k+1}\mathrm{diag}\{\frac 1{11},\frac 1{22},\frac 1{33},\frac
1{44)}\}$, $S_k=\mathrm{diag}\{\frac 1{35(k+1)},\frac 1{70(k+1)}\}$,
$S=\mathrm{diag}\{\frac 1{60},\frac 1{120}\}$. 

Note that $\mathop{\mathrm{rank}}
  \begin{smallmatrix}
    F_{2k+1}\\H_{2k+1}
  \end{smallmatrix}=2$ and $I_{2k+1}=1$ so that $N_Q=\{\ell:Q_N^+Q_N\ell=0\}$
  is nontrivial: $\ell=(0,0,1)$ belongs to $N_Q$. Theorem~\ref{t:mnmx} implies
  $\ell(\mathscr G^{2k+1}_y)=(-\infty,+\infty)$ so that the a-posteriori
  minimax error in the direction $\ell$ is infinite. Thus the estimation error
  is unbounded in general case. Really $(\ell,Q_{2k+1}^+r_{2k+1})=0$ so that $$
|(\ell,Q_{2k+1}^+r_{2k+1}-x_{2k+1})|=|x_{3,2k+1}|
$$ Note, that any function $k\mapsto x_{3,k}$ satisfies~\eqref{eq:st:ne}. In
this sense model~\eqref{eq:st:ne}-\eqref{eq:obs:ne} is non-causal. Since the
estimation error in the direction $\ell$ coincides with $x_{3,k}$ its
evolution is unpredictable for odd $k$. In this sense the subspace
$\{\alpha\ell,\alpha\in\mathbb R^1\}$ in the system state space is not
observable for odd $k$. On the other hand $I_{2k}=0$. Thus $R(Q_{2k})=\mathbb
R^3$ and the system state space is observable in any direction $\ell=\mathbb
R^3$ due to Theorem~\ref{t:mnmx}.%  Also the minimal eigen value
%   of $H'_{2k}R_{2k}H_{2k}$ grows while $k\to\infty$. Corollary~\ref{c:asymt}
%   implies that estimation error $\|x_k-\hat x_k\|$ tends to zero for even $k$

\begin{figure}[h]\centering
\begin{minipage}[c]{550pt} 
\includegraphics[viewport=1 135 700 600,width=400pt,clip,angle=270]{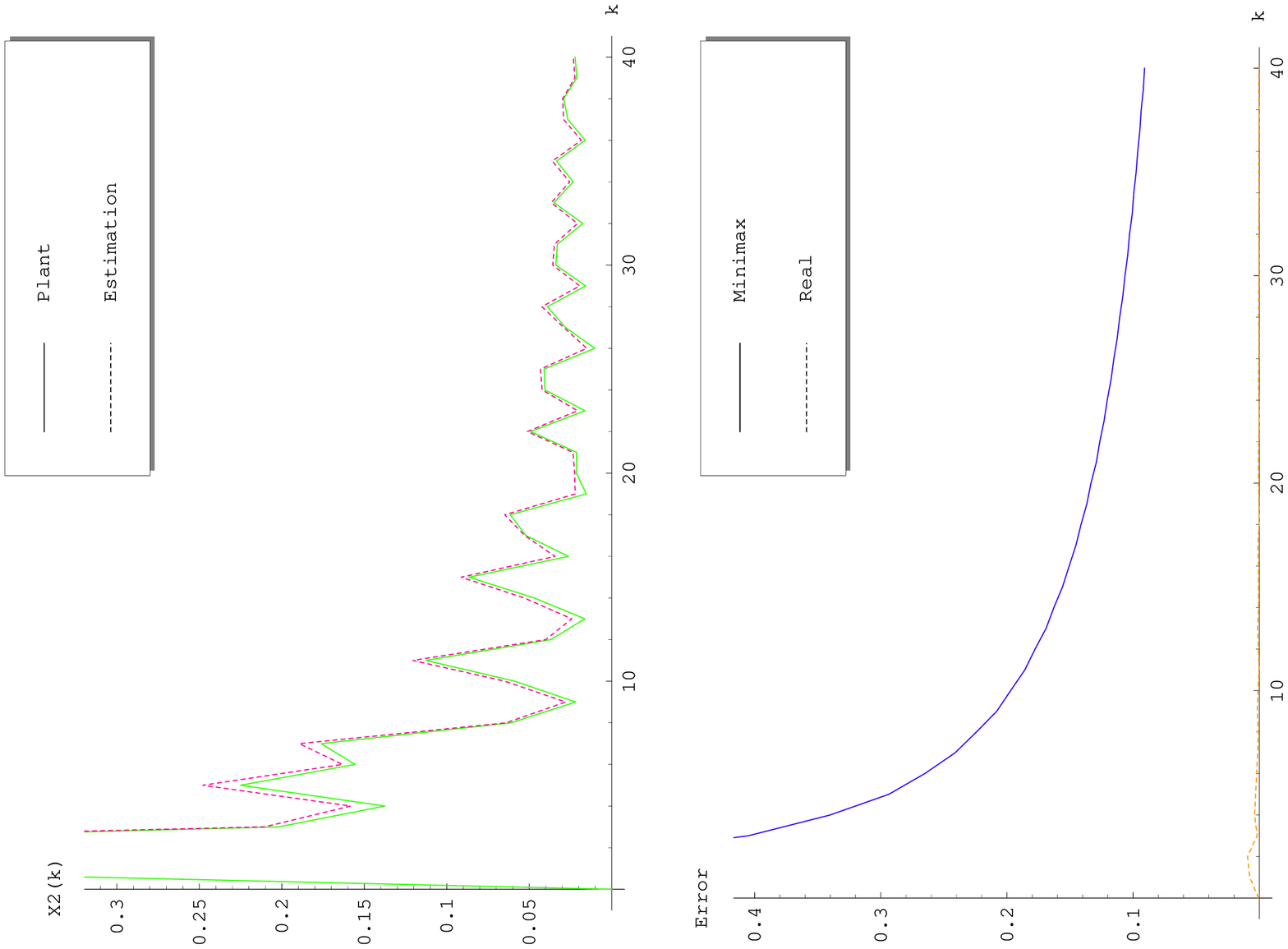}
\end{minipage}
\caption{$N=40$, state $x_{2,k}$ (solid), estimation $\hat{x}_{2,k}$ (dashed);
real error $|x_{2,k}-\hat{x}_{2,k}|$ (dashed), minimax error (solid). }
\label{fig:2}
\end{figure}

\begin{figure}[h]\centering
\begin{minipage}[c]{550pt} 
\includegraphics[viewport=1 135 700 600,width=400pt,clip,angle=270]{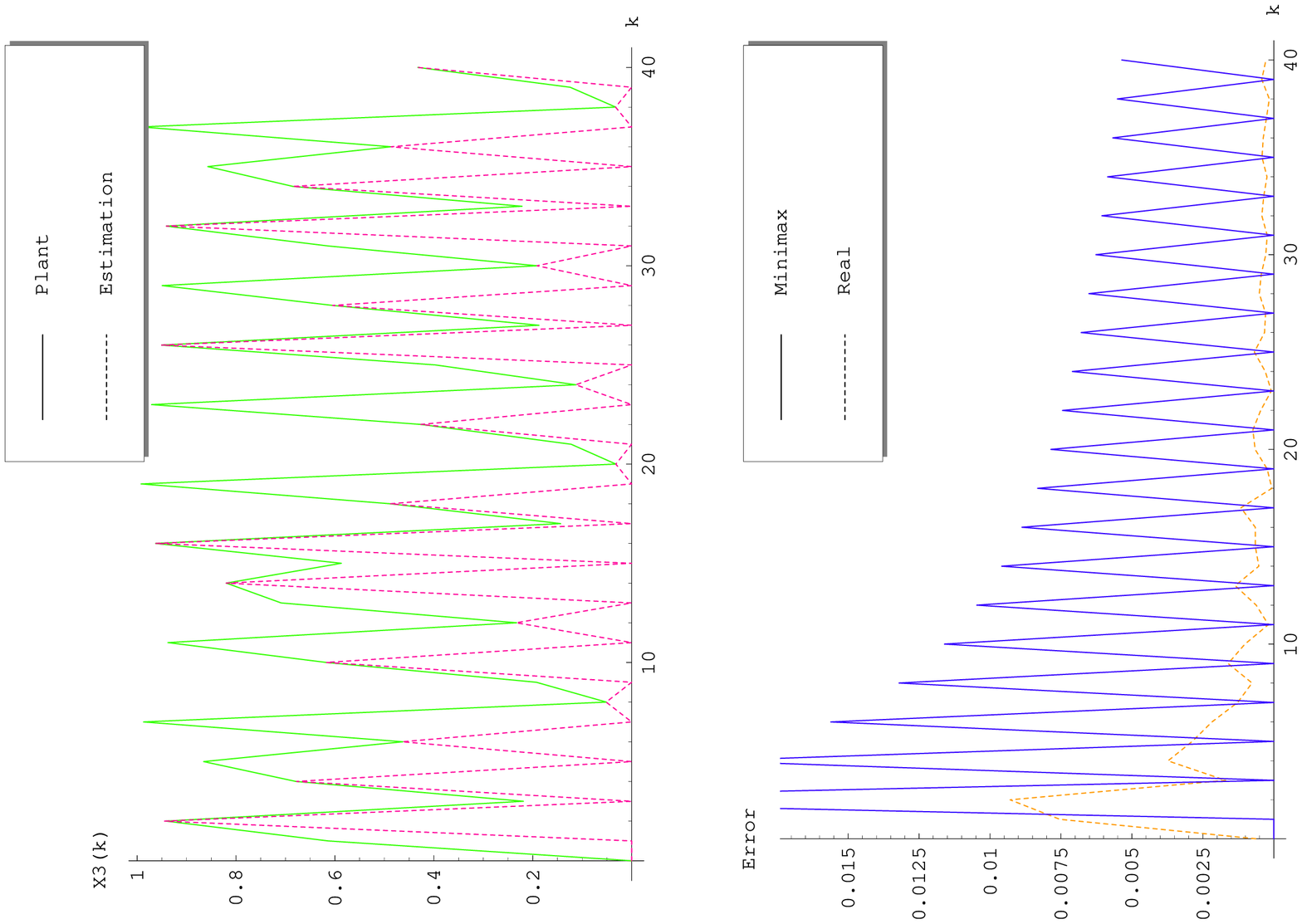}
\end{minipage}
\caption{$N=40$, state $x_{3,k}$ (solid), estimation $\hat{x}_{3,k}$ (dashed);
real error $|x_{3,k}-\hat{x}_{3,k}|$ (dashed), minimax error (solid). }
\label{fig:3}
\end{figure}

The dynamics of
$x_{i,k},\hat{x}_{i,k}$,  $|x_{i,k}-\hat{x}_{i,k}|$, $i=1,3$ and minimax
error is illustrated by figures~\ref{fig:2}-\ref{fig:3}. Note, that
Figure~\ref{fig:3} demonstrates a singular 
case: for even $k$ minimax estimation and error vanish but for odd $k$ they
gives nontrivial approximation. Thus one can observe some kind of oscillation
of the estimation curve: $\hat x_{3,2k}$ is near $x_{3,2k}$ and $\hat
x_{3,2k+1}=0$. Note, that although minimax error in the direction
$\ell=(0,0,1)$ is infinite but $(Q^+_{2k+1}\ell,\ell)=0$. Thus the
corresponding minimax error curve gives an upper bound of the
$|x_{3,k}-\hat{x}_{3,k}|$ for odd $k$ and vanishes for even $k$. 

\bibliography{refs,myrefs}

\begin{thebibliography}{10}

\bibitem{Lewis1986}
F.~Lewis.
\newblock A survey of linear singular systems.
\newblock {\em Circuits Systems Signal Process}, 5(1), 1986.

\bibitem{Mehrmann2005}
V.~Mehrmann and T.~Stykel.
\newblock Descriptor systems: a general mathematical framework for modelling,
  simulation and control.
\newblock Technical Report 292-2005, DFG Research center Matheon, Berlin, 2005.
\newblock www.matheon.de.

\bibitem{Rutkas1975}
A.~Rutkas.
\newblock Cauchy's problem for the equation $ax'(t)+bx(t)=f(t)$.
\newblock {\em Diff.Equations (translation of Diff. uravn.)}, 11:1486--1497,
  1975.

\bibitem{Dai1989}
L.~Dai.
\newblock Singular control systems.
\newblock {\em Lect.notes in Control and Information Scienses}, 8:12--24, 1989.

\bibitem{Campbell1983}
S.L. Campbell and L.R Petzold.
\newblock Canonical forms and solvable singular systems of differential
  equations.
\newblock {\em SIAM J.Alg.Discrete Methods}, 4:517--521, 1983.

\bibitem{Campbell1987}
S.~Campbell.
\newblock A general form for solvable linear time varying singular systems of
  differential equations.
\newblock {\em SIAM J. Math. Anal.}, 18(4), 1987.

\bibitem{Reich1990}
S.~Reich.
\newblock On a geometrical interpretation of differential-algebraic equations.
\newblock {\em Circuits Systems Signal Process}, 9(4), 1990.

\bibitem{Rabier1994}
P.~Rabier and W.~Rheinboldt.
\newblock A geometric treatment of implicit differential-algebraic equations.
\newblock {\em J.Differential Equtions}, 109:110--146, 1994.

\bibitem{Ilchmann2005}
A.~Ilchmann and V.~Mehrman.
\newblock A behavioral approach to time-varying linear systems. part 2:
  Descriptor systems.
\newblock {\em SIAM J. Control Optim.}, 44:1748--1765, 2005.

\bibitem{Zhuk2007}
S.~Zhuk.
\newblock Closedness and normal solvability of an operator generated by a
  degenerate linear differential equation with variable coefficients.
\newblock {\em Nonlinear Oscillations}, 10(4):1--18, 2007.

\bibitem{Ozcaldiran1989}
K.~Ozcaldiran and F.~Lewis.
\newblock Generalized reachability subspaces for singular systems.
\newblock {\em SIAM J. Control and Optimization}, 27:495--510, 1989.

\bibitem{Kurina2007}
G.~Kurina and R.~M\"arz.
\newblock Feedback solutions of optimal control problems with dae constraints.
\newblock {\em SIAM J. Control Optim.}, 46(4):1277--1298, 2007.

\bibitem{Rutkas2008}
A.~Rutkas.
\newblock Solvability of semilinear differential equations with singularity.
\newblock {\em Ukr.Math.J.}, 60, 2008.

\bibitem{Rabier1989}
P.~Rabier.
\newblock Implicit differential equations near a singular point.
\newblock {\em J. Math. Analysis Appl.}, 144:425--449, 1989.

\bibitem{Grimm1988}
J.~Grimm.
\newblock Realization and canonicity for implicit systems.
\newblock {\em SICON}, 26(6), November 1988.

\bibitem{Balakrishnan1984}
A.~Balakrishnan.
\newblock {\em Kalman Filtering Theory}.
\newblock Opt. Soft., Inc., N.Y., 1984.

\bibitem{nikh1}
R.~Nikoukhah, S.L. Campbell, and F.~Delebecque.
\newblock Kalman filtering for general discrete-time linear systems.
\newblock {\em IEEE Transactions on Automatic Control}, 44:1829--1839, 1999.

\bibitem{ishixara2}
J.Y. Ishihara, M.H. Terra, and J.C.T. Campos.
\newblock Optimal recursive estimation for discrete-time descriptor systems.
\newblock {\em Int. J. of System Science}, 36(10):1--22, 2005.

\bibitem{Deng1999}
Z.~Deng and Liu Y.
\newblock Descriptor kalman estimators.
\newblock {\em Int. J. of System Science}, 30:1205--1212, 1999.

\bibitem{Basar1995}
T.~Ba\c{s}ar and P.~Bernhard.
\newblock {\em $H_\infty$-optimal Control and Related Minimax Design Problems}.
\newblock Springer, 1995.

\bibitem{Shaked1992}
U.~Shaked and Y.~Theodor.
\newblock $h_\infty$-optimal estimation: a tutorial.
\newblock In {\em Proc. of 31st IEEE Conf. Decision and Control}, pages
  2278--2286, 1992.

\bibitem{Sayed1996}
A.~Sayed, Th. Kailath, and Hassibi B.
\newblock Linear estimation in krein spaces - part 2: Applications.
\newblock {\em IEEE Trans. on Automat. Contr.}, 41:34--50, 1996.

\bibitem{Zhang2003}
H.~Zhang, L.~Xie, and Y.~Soh.
\newblock Risk-sensitive filtering, prediction and smoothing for discrete-time
  singular systems.
\newblock {\em Automatica}, 39:57--66, 2003.

\bibitem{Xu2007}
S.~Xu and J.~Lam.
\newblock Reduced-order $h_\infty$ filtering for singular systems.
\newblock {\em System \& Control Letters}, 56(1):48--57, 2007.

\bibitem{Xu2006}
S.~Xu and J.~Lam.
\newblock Robust control and filtering of singular systems.
\newblock {\em Lect.notes in Control and Information Scienses}, 332:1--231,
  2006.

\bibitem{Bertsekas1971}
D.P. Bertsekas and I.~B. Rhodes.
\newblock Recursive state estimation with a set-membership description of the
  uncertainty.
\newblock {\em IEEE Trans. Automat. Contr.}, AC-16:117--128, 1971.

\bibitem{Nakonechnii1978}
A.~Nakonechny.
\newblock Minimax estimation of functionals defined on solution sets of
  operator equations.
\newblock {\em Arch.Math. 1, Scripta Fac. Sci. Nat. Ujer Brunensis}, 14:55--60,
  1978.

\bibitem{Tempo1985}
M.~Milanese and R.~Tempo.
\newblock Optimal algorithms theory for robust estimation and prediction.
\newblock {\em IEEE Trans. Autom. Contr.}, 30(8):730--738, Aug 1985.

\bibitem{Nakonechnii2004}
A.~Nakonechny.
\newblock Uncertain parameter estimation problems.
\newblock In {\em Scientific notes of Kyiv National University}, Cybernetics
  Faculty. KNU Press, 2004.
\newblock in ukrainian.

\bibitem{Chernousko1994}
F.~L. Chernousko.
\newblock {\em State Estimation for Dynamic Systems .}
\newblock Boca Raton, FL: CRC, 1994.

\bibitem{Milanese1991}
M.~Milanese and A.~Vicino.
\newblock Optimal estimation theory for dynamic systems with set membership
  uncertainty: An overview.
\newblock {\em Automatica}, 27:997--1009, 1991.

\bibitem{Kurzhanski1997}
A.~Kurzhanski and I.~Valyi.
\newblock {\em Ellipsoidal Calculus for Estimation and Control}.
\newblock Birkh\"auser, Boston, 1997.

\bibitem{Bakan2003}
G.~Bakan.
\newblock Analytical synthesis of guaranteed estimation algorithms of dynamic
  process states.
\newblock {\em J. of Automation and Information Sci.}, 35(5), 2003.

\bibitem{Kuntsevich1992}
V.M. Kuntsevich and M.M. Lychak.
\newblock {\em Guaranteed estimates, adaptation and robustness in control
  system}.
\newblock Springer, 1992.

\bibitem{Zhuk2008e}
S.~Zhuk.
\newblock State estimation for dynamic system described by linear operator
  equation in hilbert space.
\newblock In {\em Abs. Int. Workshop "Problems of deicision making under
  uncertainty"}, 2008.

\bibitem{Zhuk2008b}
S.~Zhuk.
\newblock State estimation for dynamical system described by uncertain linear
  equation in hilbert space.
\newblock {\em Ukrainian mathematical journal}, 61(2):178--194, 2009.

\bibitem{Zhuk2008d}
S.~Zhuk.
\newblock State estimation for dynamic system described by linear equation with
  closed operator in hilbert space.
\newblock In {\em Abs. Int. Conf. "Differential Equations and Topology"}, 2008.

\bibitem{Zhuk2009}
S.~Zhuk.
\newblock Minimax state estimation for linear equations with unbounded operator
  in hilbert space.
\newblock {\em Tavrian Bulletin of Mathematics and Informatics}, (1), 2009.
\newblock in press, Tavrian Bulletin of Mathematics and Informatics.

\bibitem{Zhuk2006c}
S.~Zhuk.
\newblock Minimax state estimation for linear descriptor systems.
\newblock Kyiv University Press, Kyiv, 2006.
\newblock Abs.of PhD thesis.

\bibitem{Zhuk2006b}
S.~Zhuk.
\newblock {\em Minimax state estimation for linear descriptor systems}.
\newblock PhD thesis, National Taras Shevchenko University of Kyiv, 2006.

\bibitem{Zhuk2004}
S.~Zhuk.
\newblock Minimax estimation of solutions of systems of linear algebraic
  equations with singular matrices.
\newblock {\em J. of Automation and Information Sci.}, 36:35--43, 2004.

\bibitem{Zhuk2004a}
S.~Zhuk.
\newblock A posteriori minimax state estimation for linear algebraic equations
  with singular matrix.
\newblock {\em Bulletin of Kyiv Univ}, 3:211--214, 2004.

\bibitem{Zhuk2003}
S.~Zhuk.
\newblock Multicriterion state estimation for linear algebraic equations with
  uncertain parameters.
\newblock In {\em Abs. Int. Conf. "Problems of decision making under
  uncertainty"}, 2003.

\bibitem{Zhuk2004b}
S.~Zhuk.
\newblock Minimax state estimation for linear singular algebraic equations.
\newblock In {\em Abs. Int. Workshop "Problems of decision making under
  uncertainty"}, pages 126--128, Ternopil, Ukraine, 2004.

\bibitem{Zhuk2004c}
S.~Zhuk.
\newblock Minimax a posteriori state estimation for singular linear algebraic
  equations.
\newblock In {\em Recent problems of modelling, forecasting and optimization}.
  Kyiv National University Press, 2004.

\bibitem{Zhuk2009a}
S.~Zhuk.
\newblock Recursive state estimation for non-causal discrete-time linear
  descriptor systems.
\newblock In {\em Proc. of IFAC Workshop Cont. Appl. Opt.}, 2009.
\newblock to appear.

\bibitem{Zhuk2009b}
S.~Zhuk.
\newblock Recursive set-membership state estimation for linear non-causal
  time-variant differential-algebraic equation with continuos time.
\newblock preprint, submited to SICON, 2009.

\bibitem{Zhuk2005c}
S.~Zhuk.
\newblock On minimax mean-square estimations for descriptor systems.
\newblock In {\em Abs. Int. Conf. "Problems of decision making under
  uncertainty"}, 2005.

\bibitem{Zhuk2005a}
S.~Zhuk.
\newblock Minimax estimation for linear descriptor differential equations.
\newblock {\em J.Appl. Comp. Math.}, 2:39--46, 2005.

\bibitem{Zhuk2006a}
S.~Zhuk.
\newblock Guaranteed estimations for linear descriptor systems.
\newblock In {\em Abs. Int. Workshop "Problems of decision making under
  uncertainty"}, 2006.

\bibitem{Zhuk2008a}
S.~Zhuk.
\newblock Minimax recursive state estimation for linear discrete-time
  descriptor systems.
\newblock in press, Int.J.System research and Informational Technologies, 2009.

\bibitem{Zhuk2006}
S.~Zhuk.
\newblock Guaranteed estimations for linear descriptor difference systems.
\newblock In {\em Abs. Int. Conf. "Problems of decision making under
  uncertainty"}, 2006.

\bibitem{Zhuk2008c}
S.~Zhuk.
\newblock Recursive state estimation for noncausal discrete-time descriptor
  systems with uncertain parameters.
\newblock In {\em Abs. Int. Conf. CDDEA}, Strechno, Slovakia, 2008.

\bibitem{Zhuk2005b}
S.~Zhuk.
\newblock Minimax state estimation for singular linear difference equations.
\newblock {\em Tavrian Bulletin of Mathematics and Informatics}, 1:16--24,
  2005.

\bibitem{Zhuk2005}
S.~Zhuk.
\newblock Minimax state estimation for linear descriptor equations with
  discrete time.
\newblock {\em Tavrian Bulletin of Mathematics and Informatics}, 2:12--22,
  2005.

\bibitem{Zhuk2008}
S.~Zhuk.
\newblock Guaranteed state estimation for linear descriptor difference
  equations.
\newblock {\em Bulletin of Kyiv Univ}, 1:102--105, 2008.

\bibitem{Zhuk2007a}
S.~Zhuk, S.~Demidenko, and A.~Nakonechniy.
\newblock Minimax state observation in linear one-dimentional 2-point boundary
  value problems.
\newblock {\em Tavrian Bulletin of Mathematics and Informatics}, 1:7--24, 2007.

\bibitem{Zhuk2008f}
S.~Zhuk, S.~Demidenko, and A.~Nakonechniy.
\newblock Minimax mean-square estimations of trends in linear regression
  problems.
\newblock In {\em Abs. Int. Conf. "Problems of decision making under
  uncertainty"}, 2008.

\bibitem{Brammer1989}
K.~Brammer and G.~Siffling.
\newblock {\em Kalman Bucy Filters}.
\newblock Artech House Inc., Norwood MA, USA, 1989.

\bibitem{Bender1987}
D.~Bender and A.~Laub.
\newblock The linear-quadratic optimal regulator for descriptor system:
  discrete-time case.
\newblock {\em Automatica}, 23:71--85, 1987.

\bibitem{Rockafellar1970}
R.~Rockafellar.
\newblock {\em Convex analysis}.
\newblock Princeton, 1970.

\bibitem{Baras1995}
J.~S. Baras and A.B. Kurzhanski.
\newblock Nonlinear filtering: The set-membership (bounding) and the $h_\infty$
  techniques.
\newblock In {\em Proc. 3rd IFAC Symp.Nonlinear Control Sys.Design}. Pergamon,
  1995.

\bibitem{Sayed2001}
A.~Sayed.
\newblock A framework for state-space estimation with uncertain models.
\newblock {\em IEEE Trans. Autom. Contr.}, 46:998--1013, 2001.

\end{thebibliography}
\bibliographystyle{unsrt}
\end{document}